\documentclass[11pt, reqno]{amsart}

 \usepackage{amsmath}
 \usepackage{amssymb}
\usepackage{bm}

 \usepackage{color}

\newcommand{\mysection}[1]{\section{#1}
      \setcounter{equation}{0}}

\newtheorem{theorem}{Theorem}[section]
\newtheorem{lemma}[theorem]{Lemma}

\theoremstyle{definition}

\theoremstyle{remark}
\newtheorem{remark}{Remark}[section]

\newcommand{\loc}{\text{\rm loc}}

 \makeatletter 
 \def\dashint{%
 \operatorname%
 {\,\,\text{\bf--}\kern-.98em\DOTSI\intop\ilimits@\!\!}}
\def\dashnorm{\,\,\text{\bf--}\kern-.5em\|}
 \makeatother

\newcommand\bbeta{\text{\raise-.2ex\hbox{$\bm{\beta}$}}}

\newcommand\osc{\operatornamewithlimits{osc}}
\newcommand\uR{\underline{R}}

\newcommand\bR{\mathbb{R}}

\newcommand\bC{\mathbb{C}}

\newcommand\cF{\mathcal{F}}

\begin{document}

\title[Heat equation  
with drift in $L_{d+1}$]
{On the heat equation  
with drift in $L_{d+1} $}

\author{N.V. Krylov}
 
\email{nkrylov@umn.edu}
\address{127 Vincent Hall, University of Minnesota,
 Minneapolis, MN, 55455}
 


\begin{abstract}
In this paper we deal with the heat equation
with drift in $L_{d+1}$. Basically, we prove that,
if the free term is in $L_{q}$ with high enough $q$,
then the equation is uniquely solvable
in a rather unusual class of functions
 such that $\partial_{t}u,
D^{2}u\in L_{p}$ with $p<d+1$ and $Du\in L_{q}$.
\end{abstract}

\maketitle

\mysection{Introduction and first main result}

Let $\bR^{d}$ be a  Euclidean space of points
$x=(x^{1},...,x^{d})$, $d\geq 2$. Define
$\bR^{d+1}=\{(t,x):t\in\bR,x\in\bR^{d}\}$ and
for $  R>0$, $(t,x)\in\bR^{d+1}$ introduce
$$
B_{R}(x)=\{y\in\bR^{d}:|y-x|<R\},\quad
B_{R}=B_{R} (0) ,\quad C_{R}=[0, R^{2})\times B_{R},
$$
$$
C _{ R}(t,x)=C _{ R}+(t,x) .
$$

Let $b(t,x)$ be Borel $\bR^{d}$-valued function
on $\bR^{d+1}$
such that for any 
  $R>0$,
  $(t,x)\in\bR^{d+1}$
\begin{equation}            \label{5.10.1}
  \|b\|^{d+1}_{L_{d+1}( C _{ R}(t,x))}
\leq \bar b_{R}^{d+1}R ,
\end{equation}
where $\bar b_{R}$, $R>0$,    
is a continuous nondecreasing 
{\em bounded\/} function.

For $f\in L_{q}(\bR^{d+1})$
vanishing for $t\geq1$ we want to investigate the equation
\begin{equation}
                                      \label{12.23.1}
\partial_{t}u+\Delta u+b^{i}D_{i}u=f
\end{equation}
in the class of functions $u\in \bigcup_{T>0}W^{1,2}_{p}((-T,1)
\times\bR^{d})$ such that $u=0$ for $t=1$, where $p<d+1$,
$q$ is large enough, $\partial_{t}=\frac{\partial}{
\partial t}$, $D_{i}=\frac{\partial}{
\partial x^{i}}$.

A somewhat unusual feature of this problem
is that $b^{i}D_{i}u\not\in L_{p}((0,1)\times\bR^{d})$
for arbitrary $u\in W^{1,2}_{p}((0,1)\times\bR^{d})$ even
vanishing for $t=1$. Therefore, if we solve 
\eqref{12.23.1} and plug the solution into
an equation with different $b$ of the same class,
we will generally not obtain a function in $L_{q}$
even locally. The author is aware of only three
similar occasions for equation with the drift term
this time growing linearly in $x$, when the solutions are sought
for in usual H\"older or Sobolev spaces without weights.
These are found in \cite{Ce_96}, \cite{KP_10}, \cite{Kr_11}.
As there, the phenomenological explanation of why
$b^{i}D_{i}u$ can be controlled is that,
as a {\em solution\/}, $u$ admits a 
probabilistic representation which
shows that, if in some direction the drift is very big,
the solution along the drift is almost constant,
so that the gradient is almost orthogonal to the drift.
This argument does not work if $u$ is just 
any arbitrary function and it shows that
$b^{i}D_{i}u$ should not be treated as a perturbation
but rather as an integral part of the operator
$L=\partial_{t}+\Delta +b^{i}D_{i}$.
This is the main reason why we concentrate
on first estimating $Du$.

Here is our first main result.
For $T\in(0,\infty)$ set $\bR^{d}_{T}=
(0,T)\times\bR^{d}$. 
\begin{theorem}
                                      \label{theorem 12.23.1}
Additionally to \eqref{5.10.1}
suppose that
$$
\|b\|_{L_{d+1}(\bR^{d}_{1})}<\infty.
$$
Let $p\in(1,d+1)$ and 
$$
q=q_{p} :=\frac{p(d+1)}{d+1-p}.
$$
Let $f$ have support in $C_{1}$ and
belong to $L_{q}(C_{1})$. Then there exists
$\hat b=\hat b(d,p)>0$ such that if $\bar b_{\infty}
\leq \hat b$, then equation \eqref{12.23.1}
has a unique solution such that
$$
\partial_{t} u,\Delta u  \in  
L_{p}( \bR^{d}_{1}),
\quad Du  \in 
L_{q}( \bR^{d}_{1}),
$$
and $u(1,\cdot)=0$. Furthermore, there exist
constants $N_{1}=N_{1}(d,p)$ and $N_{2}=N_{1}
\|b \|_{L_{d+1}( \bR^{d}_{1})}$
such that
$$
\|\partial_{t}u,D^{2}u\|_{L_{p }(\bR^{d}_{1})}
\leq N_{2}\|f\|_{L_{q}(\bR^{d}_{1})}+
 N_{1}\|f\|_{L_{p }(\bR^{d}_{1})},
$$
$$
\|Du\|_{L_{q}(\bR^{d}_{1})}
\leq N_{1}\|f\|_{L_{q}(\bR^{d}_{1})}.
$$

\end{theorem}

\begin{remark}
                                   \label{remark 12.7.1}
1. In the second part of the paper we  
relax the condition
$\bar b_{\infty}
\leq \hat b$ to $\bar b_{0+}
\leq \hat b$  and allow $f$ to be any function
in $ 
L_{q}\cap L_{p} $ but $q>q_{p}$.
The arguments there are based on some results
for diffusion processes with measurable coefficients
and are better adapted to be generalized for
fully nonlinear parabolic equations with VMO
main part.

2. From our proofs one can see that one can replace
\eqref{5.10.1} with the requirement that $b$
belongs to more general Morrey classes.
We prefer \eqref{5.10.1} for only one reason that
in the second part of the paper we use some results
from \cite{Kr_21_2} which are proved, so far,
only for $b\in L_{p,q}$ with $d/p+1/q=1$ satisfying
a condition which becomes \eqref{5.10.1} if $p=q$.

3. Hongjie Dong kindly showed the author the way to prove
the existence part in 
Theorem \ref{theorem 12.23.1} by using the theory
of parabolic Morrey's spaces. This way is probably
the one G. Lieberman had in mind writing
his Theorem  25 in \cite{Li_03}
(without proof). However, as far as the author understands,
this theorem
does not  cover Theorem \ref{theorem 12.23.1}
let alone Theorem \ref{theorem 2.16.1}
in what concerns the range of parameters.  
\end{remark}

\begin{remark}

Once we have \eqref{5.10.1} the smallness
can be always achieved by replacing $b$ with
$\lambda b$, where $\lambda $ is sufficiently small.

Also note that \eqref{5.10.1} does not imply
higher summability of $b$. For instance,
take  $\alpha\in(0,d),\beta\in(0,1)$ such that
$\alpha+2\beta=d+1$ 
and also take a continuous
bounded function $h(\tau)$, $\tau \geq 0$, with
$h(0)=0$
and consider the function
$g(t,x)=|t|^{-\beta}|x|^{-\alpha}h(|x|)$. Observe that
$$
\int_{C_{ \rho}(t,x)}g(s,y)\,dyds  
=\rho \int_{C_{1}(t',x')} 
|s|^{-\beta}|y|^{-\alpha}h(\rho|y|)\,dyds,
$$
where $t'=t/\rho^{2}$, $x'=x/\rho$. It is not hard to see that
 the last integral
is a bounded function of $(\rho,t',x')$
which tends uniformly to zero as $\rho\downarrow 0$.
 Hence, the function
$b=g^{1/(d+1)}$ satisfies \eqref{5.10.1} and even $\bar b_{0+}
=0$.
  Also clearly
for any $p>d+1$ one can find $h$, $\alpha$ and $\beta$
above such that $g^{1/(d+1)}\not\in L_{p,\loc}$.
\end{remark}

 \begin{remark}
                                   \label{remark 12.23.1}
Theorem \ref{theorem 12.23.1} is about the solvability
of the terminal value problem with zero
terminal data. Concerning nonzero data
we refer the reader to Remark \ref{remark 2.18.2}.
\end{remark}

\mysection{Auxiliary results}
                                    \label{section 12.23.1}

Set $L_{0}=\partial_{t}+\Delta$.
If $\Gamma$ is a measurable subset of $\bR^{d+1}$ and
$f$ is a function on $\Gamma$ we denote 
$$
\dashint_{\Gamma}f\,dz =\frac{1}{|\Gamma|}
\int_{\Gamma}f\,dz,
$$
where $|\Gamma|$ is the Lebesgue measure of $\Gamma$
and $z$ stands for $(t,x)$.

\begin{lemma}
                      \label{lemma 12.3.1}
Let $v\in W^{1,2}_{1}(C _{R})$
and assume that $L_{0}v=0$ in 
$C _{R}$. Then, for $\kappa\in (0, 1/4]$,
$$
 \dashint_{C _{\kappa R}}
\dashint_{C _{\kappa R}}|Dv(z_{1})- Dv(z_{2})|
\,dz_{1}dz_{2} \leq N(d)\kappa
\dashint_{C _{  R}}|Dv(z)|\,dz.
$$
\end{lemma}

 Proof. Since $Dv$ satisfies the same equation,
it suffices to prove that
\begin{equation}            \label{12.3.2}
 \dashint_{C _{\kappa R}}
\dashint_{C _{\kappa R}}| v(z_{1})-  v(z_{2})|
\,dz_{1}dz_{2} \leq N\kappa
\dashint_{C _{  R}}| v(z)|\,dz.
\end{equation}Self-similar transformations
allows us to assume that $R=1$.

We know (see, for instance, theorem 8.4.4 of \cite{Kr_96})
that
$$
 \dashint_{C _{\kappa }}
\dashint_{C _{\kappa }}|v(z_{1})- v(z_{2})|\,dz_{1}dz_{2}
\leq N\kappa\sup_{C _{\kappa }}(|\partial_{t}v
+|Dv|)\leq N\kappa\sup_{C _{2\kappa }}|v|,
$$
where the last supremum is easily estimated through
$$
\int_{C_{1}}|v|\,dz .
$$
The lemma is proved. 

For $\pi\in(1,d+2)$
introduce
$$
\pi^{*}=\frac{\pi(d+2)}{d+2-\pi},
$$
and observe that, if $\pi<d+1$,
$$
 \pi^{*}
<\pi(d+1)/(d+1-\pi)=:q_{\pi} ,
$$ 
$\|b^{i}D_{i}u\|_{L_{\pi}}\leq \|b\|_{L_{d+1}}  
\|Du\|_{L_{q_{\pi}}}$, whereas by embedding
 theorems $\partial_{t}u,D^{2}u
\in L_{\pi}$ only implies that $Du\in L_{\pi^{*}}$.
This
presents the main obstacle on the way of
``usual'' Sobolev space PDE theory for the operator $L$
when lower-order terms are treated as
perturbations.

Define $\partial'C_{R}=\bar C_{R}\setminus(\{t=0\}
\times B_{R})$ and introduce the notation
$$
\dashnorm g \|^{r}_{L_{r}(C_{R})}
=\dashint_{C_{R}}|g|^{r}\,dz.
$$

\begin{lemma}
                      \label{lemma 12.3.10}
Let   $w\in W^{1,2}_{\pi}(C _{R})$
and assume that $L_{0}w=f$ in 
$C _{R}$ and $w=0$ on $\partial'C_{R}$. Then
\begin{equation}            \label{12.5.10}
 \dashnorm Dw \|_{L_{\pi^{*}}(C_{R})}\leq N(d,\pi)R
\dashnorm f \|_{L_{\pi}(C_{R})}.
\end{equation}
\end{lemma}

Proof. Rescailing allows us to assume that $R=1$.
In that case the $W^{1,2}_{\pi}(C_{1})$-norm of $w$
is estimate through the $L_{\pi}(C_{1})$-norm of $f$.
After that it only remains to use embedding theorems.
The lemma is proved.
 
 This result is used below with $1$ in place of $\pi^*$.

\begin{lemma}
                                      \label{lemma 12.5.3}
Let   $u\in W^{1,2}_{\pi}(C _{R})$.
Introduce $L_{0}u=f$. 
Then, for $\kappa\in (0, 1/4]$, with $N=N(d,\pi)$,
$$
\dashint_{C _{\kappa R}}
\dashint_{C _{\kappa R}}|Du(z_{1})-Du(z_{2})|
\,dz_{1}dz_{2}   \leq N\kappa
\dashint_{C _{  R}}|Du(z)|\,dz
$$
\begin{equation}            \label{12.5.1}
 +
N\kappa^{-2d-4}R\Big(\dashint_{C _{  R}}|f|^{\pi}
\,dz\Big)^{1/\pi}.
\end{equation}
\end{lemma}

Proof. Introduce $v\in W^{1,2}_{\pi}(C_{R})$
such that $L_{0}v=0$ and $ v=u  $ on $\partial' C_{R}$
and let $w=u-v$.
Then $L_{0}w=L_{0}u=f $ and 
$$
\dashint_{C _{\kappa R}}
\dashint_{C _{\kappa R}}|Dv(z_{1})- Dv(z_{2})|
\,dz_{1}dz_{2} \leq N\kappa
\dashint_{C _{  R}}|Dv|\,dz,
$$
$$
\dashint_{C _{  R}}|Dv|\,dz\leq
\dashint_{C _{  R}}|Du|\,dz+
\dashint_{C _{  R}}|Dw|\,dz 
$$
$$
\leq \dashint_{C _{  R}}|Du|\,dz
+NR
\Big(\dashint_{C _{  R}}|f|^{\pi}\,dz\Big)^{1/\pi},
$$
$$
\dashint_{C _{\kappa R}}
\dashint_{C _{\kappa R}}|Dw(z_{1})- Dw(s_{2})|\,dz_{1}dz_{2}
 \leq N\kappa^{-2d-4}R
\Big(\dashint_{C _{  R}}|f|^{\pi}\,dz\Big)^{1/\pi}.
$$
 These  computations imply \eqref{12.5.1}
and the lemma is proved.

\begin{theorem}
                                      \label{theorem 12.5.1}
Let $\pi\in(1,d+1)$ and $u\in W^{1,2}_{\pi}(C _{R})$.
Set $f=Lu$. 
Then, for $\kappa\in (0, 1/4]$, with $N=N(d,\pi)$,
$$
\dashint_{C _{\kappa R}}
\dashint_{C _{\kappa R}}|Du(z_{1})-Du(z_{2})|\,dz_{1}dz_{2} 
\leq N\kappa
\dashint_{C _{  R}}|Du(z)|\,dz
$$
\begin{equation}            \label{12.5.01}
 +
N\bar b_{R}\kappa^{-2d-4}\Big(\dashint_{C _{  R}}|Du|^{q_{\pi} }
\,dz\Big)^{1/q_{\pi} }+
NR\kappa^{-2d-4}\Big(\dashint_{C _{  R}}
|f|^{\pi}\,dz\Big)^{1/\pi}.
\end{equation}
\end{theorem}

This result follows from \eqref{12.5.1} and the fact that
by H\"older's inequality
$$
\Big(\dashint_{C _{  R}}|b|^{\pi}\,|Du|^{\pi}\,dz\Big)^{1/\pi}
$$
$$
\leq \Big(\dashint_{C _{  R}}|b|^{d+1}\,dz\Big)^{1/(d+1)}
\Big(\dashint_{C _{  R}}|Du|^{q_{\pi} }\,dz
\Big)^{1/q_{\pi} }.
$$  

The last term in \eqref{12.5.01} presents
certain inconvenience which forced us to
assume that $f=0$ outside $C_{1}$.

\begin{lemma}
                                    \label{lemma 12.12.3}
Let $g\geq0$ have support in $C_{1}$
and be integrable. Let  $z \in \bR^{d+1}$,
$\kappa\in (0, 1/4]$. Then
for any $R>0$ and $z_{0}\in C_{\kappa R}(z)$
$$
 R^{\pi}\dashint_{C_{R}(z)}g\,dxdt
\leq N Mg(z_{0})+N(|z_{0}|+1)^{\pi-d-2}\int_{C_{1}}g\,dxdt,
$$
where $Mg$ is the parabolic Hardy-Littlewood
maximal function,
$|z_{0}|=\sqrt{|t_{0}|}+|x_{0}|$, and $N=N(d,\pi)$.
\end{lemma}

Proof. Introduce $\hat C:=C_{2}(-1,0)$ which is a cylinder
strictly containing $C_{1}$ and consider a few cases.

{\em Case $z_{0}\in \hat C$\/}. If $R\leq 1$, then by definition
$$
 R^{\pi}\dashint_{C_{R}(z)}g\,dxdt\leq Mg(z_{0}).
$$
However, if $R>1$, then
$$
R^{p}\dashint_{C_{R}(z)}g\,dxdt\leq NR^{\pi-d-2}
\int_{C_{1}}g\,dxdt\leq N\int_{C_{1}}g\,dxdt\leq N
Mg(z_{0}).
$$

{\em Case $z_{0}\not\in \hat C, t_{0}\leq-1$\/}.
In this case in order for the intersection
of $C_{R}(z)$ and $C_{1}$ be nonempty we have to have
$t_{0}+R^{2}>0$ and $|x_{0}|-2R<1$,
that is $R\geq \max(\sqrt{|t_{0}|},(1/2)(|x_{0}|-1))$.
By taking into account that $|t_{0}|\geq 1$
it is not hard to see that
$$
 \max(\sqrt{|t_{0}|},(1/2)(|x_{0}|-1))
\geq \nu\big(\sqrt{|t_{0}|}+|x_{0}|+1\big),
$$
where $\nu>0$ is an absolute constant. In that case
\begin{equation}
                                             \label{12.13.1}
R^{\pi}\dashint_{C_{R}(z)}g\,dxdt\leq NR^{\pi-d-2}
\int_{C_{1}}g\,dxdt
\leq N\frac{1}{(1+|z|)^{d+2-\pi}}
\int_{C_{1}}g\,dxdt.
\end{equation}

{\em Case $z_{0}\not\in \hat C, t_{0}\geq 3$\/}.
This time $C_{R}(z)\cap C_{1}\ne\emptyset$ only if
$1+R^{2}>t_{0}$ and $|x_{0}|-2R<1$, that is
$R\geq \max(\sqrt{ t_{0}-1},(1/2)(|x_{0}|-1))$,
which leads to \eqref{12.13.1} again.

{\em Case $z_{0}\not\in \hat C, t_{0}\in [-1, 3]$\/}.
Here $|x_{0}|\geq 2$ and $C_{R}(z)\cap C_{1}\ne\emptyset$
only if $|x_{0}|-2R<1$, that is $R\geq (1/2)(|x_{0}-1)
\geq(1/8)(|x_{0}|+1)$,
which leads to \eqref{12.13.1} again.
The lemma is proved.  

Here is the main a priori estimate.
Recall that $p\in(1,d+1)$ and $q=p(d+1)/(d+1-p)$.

\begin{lemma}
                                        \label{lemma 12.13.1}
Let  $u\in \bigcup_{T>0} 
W^{1,2}_{p}((-T,1)\times\bR^{d})$
and $Du\in  L_{q}((-\infty,1)\times\bR^{d})$. Assume that
$u(1,\cdot)=0$, 
$f:=Lu\in L_{q}((-\infty,1)\times\bR^{d})$,
and $f$ has support in $C_{1}$. Then there exists
a constant $\hat b=\hat b(d,p)>0$ such that,
if $\bar b_{\infty}\leq \hat b$, then
\begin{equation}
                                        \label{12.23.3}
\|Du\|_{L_{q}((-\infty,1)\times\bR^{d})}
\leq N\|f\|_{L_{q}(C_{1})},
\end{equation}
where $N=N(d,p )$.

\end{lemma}

Proof. We extend $u$ and $f$ as zero for $t
>1$. 
Let $\bC$ be the collection of $C_{R}(t,x)$,
$R>0$, $(t,x)\in \bR^{d+1}$.
For functions $h=h(z)$ on $\bR^{d+1}$ for which it makes sense
introduce
$$
h^{\sharp}(z)=\sup_{\substack{C\in\bC,\\
C\ni z}}
\dashint_{C }\dashint_{C }
\big|h(z_{1})-h(z_{2})\big| \,dz_{1}dz_{2} .
$$

Observe that if $z\in\bR^{d+1}$ and $z\in C\in\bC$,
then owing to Theorem \ref{theorem 12.5.1}
and Lemma \ref{lemma 12.12.3} with $\pi=(1+p)/2$
$$
\dashint_{C }\dashint_{C }
\big|Du(z_{1})-Du(z_{2})\big| \,dz_{1}dz_{2}
\leq N\kappa M|Du|(z)
$$
$$
+N\bar b_{\infty}\kappa^{-2d-4}
\Big(M\big(|Du|^{q_{\pi}}\big)(z)\Big)^{1/q_{\pi}}
+N\kappa^{-2d-4}\Big(M\big(|f|^{\pi}\big)(z)\Big)^{1/\pi}
$$
$$
+N\kappa^{-2d-4}\|f\|_{L_{\pi}(C_{1})}h(z),
$$
where $h(z)=(|z |+1)^{1-(d+2)/\pi}$. Due to the arbitrariness
of $C\ni z$ one can replace here the left-hand side
with $(Du)^{\sharp}(z)$.
Observe that,  $ \nu:=q_{\pi}((d+2)/\pi-1)
=(d+2-\pi )(d+1)/(d+1-\pi)>d+2$
and
$$
\int_{\bR^{d+1}}(|z |+1)^{q_{\pi}(1-(d+2)/\pi)}\,dz 
=N\int_{\bR^{d}}(|x|+1)^{2-\nu}\,dx<\infty.
$$

 Then by the Fefferman-Stein
theorem and by the Hardy-Littlewood maximal function
theorem (observe that $q>q_{\pi}$) we get
$$
\|Du\|_{L_{q}((-\infty,1)\times\bR^{d})}
\leq N_{1}(\kappa+\bar b_{\infty}\kappa^{-2d-4})
 \|Du\|_{L_{q}((-\infty,1)\times\bR^{d})}
$$
$$
+N\kappa^{-2d-4}\|f\|_{L_{q}(C_{1})}.
$$
To obtain \eqref{12.23.3} now it only remains
to choose first small $\kappa$ and then $\hat  b$
so that $N_{1}(\kappa+\hat  b\kappa^{-2d-4})
\leq 1/2$. The lemma  is proved.

 {\bf Proof of uniqueness in Theorem
\ref{theorem 12.23.1}}. Let $f=0$, our goal is to
show that the only solution $u$ with the specified
properties is zero. Since $L_{0}u=-b^{i}D_{i}u\in
L_{p}(\bR^{d}_{1})$, we have that
$u\in W^{1,2}_{p}(\bR^{d}_{1})$.  

Now
fix a $t_{0}>0$
close to zero, such that 
$u(t_{0},\cdot)\in W^{2}_{p}(\bR^{d})$ and 
  for $t\leq t_{0}$ define $w$
as a solution given by means of the heat semigroup
of the equation $L_{0}w=0$, $t\leq t_{0}$,
with terminal data $w(t_{0},\cdot)=u(t_{0},\cdot)$.
For $t\in[t_{0},1]$ set $w=u$.
Then $w$  is of class
$\bigcup_{T>0} 
W^{1,2}_{p}((-T,1)\times\bR^{d})$ and satisfies
$L_{0}w+I_{t>t_{0}}b^{i}D_{i}w=0$ in $(-\infty,1)
\times \bR^{d}$ with zero terminal condition.
By using the explicit representation of 
$ w $
for $t\leq t_{0}$ and the fact that
by assumption $Du\in L_{q}(\bR^{d}_{1})$,
 one easily shows that
$Dw\in L_{q}((-\infty,1)
\times \bR^{d})$. But then owing to \eqref{12.23.3},
$Dw=0$ and $L_{0}w=0$ in $(-\infty,1)
\times \bR^{d}$ and $L_{0}u=0$ in 
$(t_{0},1)
\times \bR^{d}$. It follows that $u=0$
for $t\in[t_{0},1]$ and since $t_{0}$ can be chosen
arbitrarily   close to 
 0 , $u=0$ in $\bR^{d}_{1}$,
and the uniqueness of solutions is established.

Now comes the last step needed to prove the existence
part in Theorem \ref{theorem 12.23.1}.

\begin{lemma}
                                       \label{lemma 12.28.1}
Let $f\in C^{\infty}_{0}(C_{1})$, $f=0$ outside $C_{1}$,
 and 
$b \in C^{\infty}_{0}(\bR^{d+1})$. Define $u$ as the classical
solution of $Lu=f$ for $t\leq1$ with terminal condition $u(1,
\cdot)=0$. Assume that
  $\bar b_{\infty}\leq\hat b$.
  Then   
\begin{equation}
                                        \label{12.28.3}
\|Du\|_{L_{q}(\bR^{d}_{1})}
\leq N\|f\|_{L_{q}(\bR^{d}_{1})},
\end{equation}
where $N=N(d,p )$. Furthermore,  
\begin{equation}
                                        \label{12.28.4}
\|\partial_{t}u,D^{2}u\|_{L_{p }(\bR^{d}_{1})}
\leq N_{1}\|f\|_{L_{q}(\bR^{d}_{1})}+
 N_{2}\|f\|_{L_{p }(\bR^{d}_{1})},
\end{equation}
where $N_{1}=N(d,p )\|b\|_{L_{d+1}(\bR^{d}_{1})}$,
$N_{2}=N_{2}(d,p )$.
\end{lemma}

Proof.  
The existence of smooth bounded $u$
is a classical result.  For $t\leq 0$, define $u(t,x)$ as the solution of
$L_{0}u=0$ with termunal data $u(0,\cdot)$. For $t<0$, $u(t,x)$ is
just a caloric function and it  is represented
by means of the fundamental solution of the heat equation.
 Furthermore, we have
$q>(d+2)/(d+1)$ so that simple estimates
show that $Du\in  L_{q}((-\infty,1)\times\bR^{d})$.
Now   \eqref{12.28.3} follows from Lemma \ref{lemma 12.13.1}.

Estimate \eqref{12.28.3} and H\"older's inequality show that
$$
\|b^{i}D_{i}u\|_{L_{p }(\bR^{d}_{1})}
\leq \|b\|_{L_{d+1}(\bR^{d}_{1})}\|Du\|
_{L_{q}(\bR^{d}_{1})},
$$
which implies that $f-b^{i}D_{i}u\in L_{p }(\bR^{d}_{1})$,
so that \eqref{12.28.4} is a classical result.
The lemma is proved.

{\bf Proof of Theorem \ref{theorem 12.23.1}}.
The uniqueness part is taken care of above.
To prove the existence, take $f_{n}\in C^{\infty}_{0}
(C_{1})$ converging to $f\in L_{q}(C_{1})$
and $b_{n} \in C^{\infty}_{0}(\bR^{d+1})$
converging to $b$ in $L_{d+1}(\bR^{d}_{1})$
and having $\bar b_{R}$ the same for all $n$
(just use mollifiers
and cut-off's). Then by Lemma \ref{lemma 12.28.1}
we have solutions $u_{n}$ of $L_{0}u_{n}+b^{i}_{n}
D_{i}u_{n}=f_{n}$ admitting estimates
\eqref{12.28.3} and \eqref{12.28.4} with $u_{n}$ and $f_{n}$
in place of $u$ and $f$ and with the constants
independent of $n$. Now to prove the theorem
it only remains to check that, if $Du_{n}\to Du$
weakly in $L_{q}(\bR^{d}_{1})$, then
$$
b^{i}D_{i}u^{n}\to b^{i}D_{i}u
$$
weakly in $L_{p}(\bR^{d}_{1})$. As we have seen
a few times the sequence $b^{i}D_{i}u^{n}$ is bounded in  
$L_{p}(\bR^{d}_{1})$, so we need
$$
\int_{\bR^{d}_{1}}\phi b^{i}D_{i}u^{n}\,dz
\to \int_{\bR^{d}_{1}}\phi b^{i}D_{i}u\,dz
$$
for any $\phi\in L_{p/(p-1)}(\bR^{d}_{1})$.
The latter holds indeed, since by H\"older's inequality
$\phi b\in L_{q/(q-1)}(\bR^{d}_{1})$.
The theorem is proved.

\mysection{Case when $\bar b_{R}$ is small}

We suppose that assumption \eqref{5.10.1}
is satisfied and $\|b\|_{L_{d+1}(R^{d+1})}<\infty$.
For $\delta\in(0,1)$ take the finite continuous function
$\bar N(d,d+1,\delta)$ introduced in Theorem 2.3
of \cite{Kr_21_1} and assume that there exists
$\uR\in(0,\infty)$ such that
$$
\bar N(d,d+1,1/2)\bar b_{\uR}<1.
$$
Next, let $d_{0}=d_{0}(d,1/2,  \uR )\in(d/2,d)$
be taken from \cite{Kr_21_1}. 

Below, in Theorem \ref{theorem 2.16.1} (for  $\delta=1/2$)
  $p\in [d_{0}+1,d+1)$ and
$$
q>q_{p}=\frac{p(d+1)}{d+1-p}.
$$

\begin{theorem}
                                     \label{theorem 2.16.1}
There  is a constant $\hat b>0$,  
 depending only
on $d$, $\uR$, $p$, $q$, $\bar b_{\uR}$, 
$\|b\|_{L_{d+1}(R^{d+1})}$, and 
the function $\bar N(d,d+1,\cdot)$,
 such that if 
\begin{equation}
                                               \label{2.18.1}
\bar b_{R_{0}}\leq \hat b
\end{equation}
for an $R_{0}\in(0,\uR]$, then there exists
a constant $N_{0}$, depending only on
what $\hat b$ depends on and on $R_{0}$ and $\bar b_{\infty}$,
such that
for any
$\lambda>N_{0}$ and $f\in L_{p}(\bR^{d+1})
\cap L_{q}(\bR^{d+1})$ there exists a unique
solution of $Lu-\lambda u=f$ in the class of functions
such that
\begin{equation}
                                            \label{12.11.4}
\partial_{t}u,D^{2}u\in L_{p}(\bR^{d+1}),\quad
Du\in L_{q}(\bR^{d+1}),\quad u\in L_{p}(\bR^{d+1})
\cap L_{q}(\bR^{d+1}).
\end{equation}
\end{theorem} 

We prove Theorem \ref{theorem 2.16.1}
after some preparations.
For  $\gamma\in(0,1)$  and   $\rho>0$
introduce the restricted sharp function of $h$   by
the formula   
\begin{equation}           \label{8,13.1}
 h^{\sharp}_{\gamma,\rho}(z)=\sup\big\{I_{r}( h,z_{0}):
z_{0}\in
(0,\infty)\times\bR^{d},r\in(0,\rho],C_{r}(z_{0})\ni z\big\},
\end{equation}
where
$$
I_{r}(h,z)=
\bigg(
\dashint_{C_{r}(z )}\dashint_{C_{r}(z  )}
\big|h(z_{1})-h(z_{2})\big|^{\gamma}\,dz_{1}dz_{2}\bigg)^{1/\gamma}.
$$

Here is Theorem C.2.4 of \cite{Kr_18}.
\begin{theorem}      \label{theorem 8,4.1}
Let $q\in(1,\infty)$, $\kappa\in(0,1]$, $R\in(0,\infty)$,
 and $h\in L_{q}(C _{R(1+ 2\kappa)})$.
Then   
\begin{equation}       \label{8,4.1}
\dashnorm h\|_{L_{q}(C _{R})}\leq N 
 \dashnorm  h^{\sharp}_{ \gamma,\kappa R}
 \|_{L_{q}(C _{R })}+ N\kappa^{-\chi }
 \dashnorm   h  \|_{L_{\gamma}(C _{ R })} ,
\end{equation}
where  $\chi =(d+2)/\gamma$ 
 and the constants
 $N$ depend only on $d,\gamma$, and $q$.
\end{theorem}

We also need a very particular case of Theorem 5.3
of \cite{Kr_21_2}.  
\begin{theorem}
                                  \label{theorem 12.7.1}
There is $\gamma\in(0,1)$ depending only on
$d,\uR$
such that for any $R\in(0,\uR]$, $u\in W^{1,2}_{d_{0}+1}
(C_{R})$
$$
 \dashnorm   Du  \|_{L_{\gamma}(C _{ R })} 
\leq NR\dashnorm f\|_{L_{d_{0}+1}(C_{R})}
+NR^{-1}\osc_{\partial'C_{R}} u ,
$$
where $f=Lu$ and the constants $N$ depend only on 
$ d,d_{0}, \uR $, $\bar b_{\uR}$  and
the function  $\bar N(d,d+1,\cdot)$.
 \end{theorem}

By combining this with embedding theorems
and taking into account that $d_{0}+1>d/2+1$
we come to  the following.
\begin{lemma}
                                          \label{lemma 12.8.2}
For $\gamma$ from Theorem \ref{theorem 12.7.1}
and the same type of constants $N$,
  for any $R\in(0,\uR]$ and $u\in W^{1,2}_{d_{0}+1}
(C_{R})$ we have
$$
 \dashnorm   Du  \|_{L_{\gamma}(C _{ R })} 
\leq NR\dashnorm f\|_{L_{d_{0}+1}(C_{R})}
$$
$$
+NR  \dashnorm \partial_{t}u,D^{2}u\|_{L_{d_{0}+1}(C_{R})}
+NR^{-1} \dashnorm u\|_{L_{d_{0}+1}(C_{R})},
$$
where $f=Lu$.

\end{lemma}

\begin{remark}
                                  \label{remark 12.9.1}
Below we use the fact that by H\"older's inequality
if $q\geq d_{0}+1$ and $\kappa\in(0,1]$ that
$$
\dashnorm f\|_{L_{d_{0}+1}(C_{R})}
\leq  \dashnorm f\|_{L_{q}(C_{R})}
\leq N(d)\dashnorm f\|_{L_{q}(C_{R+2\kappa R})}.
$$
\end{remark}

\begin{lemma}
                              \label{lemma 12.7.1}
Let $p\in [d_{0}+1,d+1)$ and
$$
q>q_{p}.
$$
Take $\kappa\in(0,1]$, $R\in(0,\uR]$, and
$u\in W^{1,2}_{p}(C_{R+2\kappa R})$. Set $f=Lu$. Then
$$
\dashnorm Du \|
_{L_{q}(C _{R })}\leq N(\kappa+\bar b_{R})
\dashnorm Du\|_{L_{q}(C _{R+2\kappa R})}+NR(1+\kappa^{-\chi})
\dashnorm f\|_{L_{q}(C _{R+2\kappa R})}
$$
\begin{equation}
                                          \label{12.8.1}
+N\kappa^{-\chi}R  \dashnorm \partial_{t}u,
D^{2}u\|_{L_{d_{0}+1}(C_{R})}
+N\kappa^{-\chi}R^{-1} \dashnorm u\|_{L_{d_{0}+1}(C_{R})},
\end{equation}
where the constants $N$ depend only on 
$ d,d_{0}, \uR $, $p,q$,  $\bar b_{\uR}$ and
the function  $\bar N(d,d+1,\cdot)$.
\end{lemma}

Proof.
Let $h=Du$. Then  for
$z\in C_{R}$, $r\leq R$, 
$z_{0}\in
(0,\infty)\times\bR^{d}$, and   $C_{\kappa r}(z_{0})  
\ni z$ we have $C_{ r}(z_{0})
\subset C_{R+2\kappa R}$. It follows from Theorem 
\ref{theorem 12.5.1} that
$$
I_{\kappa r}(h,z_{0})\leq N(\kappa+\bar b_{r})
\Big(\dashint_{C _{  r}(z_{0})}|h|^{q_{p}}
\,dxdt\Big)^{1/q_{p}}
$$
$$
+
NR\Big(\dashint_{C _{  r}(z_{0})}|f|^{p}\,dxdt\Big)^{1/p}.
$$
Hence,
  on $C_{R}$
$$
h^{\sharp}_{ \gamma,\kappa R}(z)
\leq N(\kappa+\bar b_{R})
\Big(\dashint_{C _{  r}(z_{0})}I_{C_{R+2\kappa R}}
|h|^{q_{p}}
\,dxdt\Big)^{1/q_{p}}
$$
$$
+
NR\Big(\dashint_{C _{  r}(z_{0})}I_{C_{R+2\kappa R}}|f|^{p}\,dxdt\Big)^{1/p}
$$
$$
\leq N(\kappa+\bar b_{R})\Big(M\Big(I_{C_{R+2\kappa R}}|h|^{
 q_{p}}
\Big)(z)\Big)^{1/q_{p}}
$$
$$
+NR\Big(M\Big(I_{C_{R+2\kappa R}}|f|^{p}
\Big)(z)\Big)^{1/p }.
$$

For $q> q_{p}$
by Hardy-Littlewood
$$
\big\|h^{\sharp}_{ \gamma,\kappa R}\big\|
_{L_{q}(C _{R })}\leq N(\kappa+\bar b_{R})
\|h\|_{L_{q}(C _{R+2\kappa R})}+NR
\|f\|_{L_{q}(C _{R+2\kappa R})} 
$$
and this along with Theorem \ref{theorem 8,4.1},
  Lemma \ref{lemma 12.8.2}, and Remark 
\ref{remark 12.9.1} yields the desired result.

Now we are going to replace $C$ with $C(z)$ in
\eqref{12.8.1} thus obtaining an inequality between two
functions on $\bR^{d+1}$ and then take the $L_{q}$-norms
of both sides as a functions on $\bR^{d+1}$.
We need a lemma.

\begin{lemma}
                                         \label{lemma 12.8.3}
Let $h$ be a nonnegative function on
$\bR^{d+1}$ and let $\infty>q\geq t\geq1$,
$r,s\in[1,\infty)$ be such that
\begin{equation}
                                          \label{12.9.2}
1+\frac{t}{q}=\frac{1}{r}+\frac{1}{s},
\end{equation}
 
Then for any $R\in(0,\infty)$
\begin{equation}
                                          \label{12.8.6}
\Big(\int_{\bR^{d+1}}\dashnorm h\|^{q}_{L_{t}(C_{R}(z))}
\,dz\Big)^{1/q}\leq
N(d)R^{-(d+2)(1-1/s)/t}
\|h^{t}\|^{1/t}_{L_{r}(\bR^{d+1})}
\end{equation}

\end{lemma}

Proof. Observe that
$$
\dashnorm h\|^{t}_{L_{t}(C_{R}(z))}
=NR^{-d-2}h^{t}*I_{C_{R}}(z).
$$
Therefore, the left-hand side of \eqref{12.8.6} is
$$
NR^{-(d+2)/t}\|h^{t}*I_{C_{R}}\|^{1/t}_{L_{q/t}(\bR^{d+1})}.
$$
By Young's inequality the $L_{q/p}$-norm of the above
convolution is dominated by
$$
\|h^{t}\|_{L_{r}(\bR^{d+1})}\|I_{C_{R}}\|_{L_{s}(\bR^{d+1})}
=NR^{(d+2)/s}\|h^{t}\|_{L_{r}(\bR^{d+1})}
$$
and the result follows.

Under the conditions of Lemma \ref{lemma 12.7.1}
 we see that \eqref{12.8.1} with $C(z)$ in place of
$C$ yields
$$
 \|Du\|_{L_{q}(\bR^{d+1})}\leq
 N(\kappa+\bar b_{R})  
\|Du\|_{L_{q}(\bR^{d+1})}
+N 
 R  \| f\|_{L_{q}(\bR^{d+1})}+I,
$$
where $I$ is the sum of the $L_{q}(\bR^{d+1})$-norms
of the last two terms in \eqref{12.8.1} with $C(z)$ in place of
$C$. To estimate these we use Lemma \ref{lemma 12.8.3}
by taking $t=d_{0}+1$, $r=p/t$, and $s>1$ defined from
\eqref{12.9.2}.  
Then we see that
$$
\Big(\int_{\bR^{d+1}}\dashnorm \partial_{t}u,
D^{2}u\|^{q}_{L_{d_{0}+1}(C_{R}(z))}\,dz\Big)^{1/q}
\leq NR^{(d+2)(1/p-1/q)}
\|\partial_{t}u,
D^{2}u\| _{L_{p}(\bR^{d+1})}.
$$
Similarly we treat the last term in \eqref{12.8.1}
and conclude that
$$
 \|Du\|_{L_{q}(\bR^{d+1})}\leq
 N_{1}(\kappa+\bar b_{R})  
\|Du\|_{L_{q}(\bR^{d+1})}
+N 
 R  \| f\|_{L_{q}(\bR^{d+1})} 
$$
$$
+N\kappa^{-\chi}
R^{1+(d+2)(1/p-1/q)}
\|\partial_{t}u,
D^{2}u\| _{L_{p}(\bR^{d+1})}
$$
\begin{equation}
                                      \label{12.9.5}
+N\kappa^{-\chi}R^{ (d+2)(1/p-1/q)-1}\|u\|_{L_{p}(\bR^{d+1})}.
\end{equation}

We fix $\kappa$ so that $N_{1}\kappa\leq 1/4$,
observe that $N_{1}$
 depends only on 
$ d,d_{0}, \uR $, $p,q$,  $\bar b_{\uR}$, and
the function  $\bar N(d,d+1,\cdot)$,
and in the future will only concentrate on $R$
such that 
$$
N_{1}\bar b_{R}\leq 1/4.
$$
 In that case
provided that the left-hand side of \eqref{12.9.5}
is finite we get

$$
 \|Du\|_{L_{q}(\bR^{d+1})}\leq
 N 
 R  \| f\|_{L_{q}(\bR^{d+1})} 
+N 
R^{1+(d+2)(1/p-1/q)}
\|\partial_{t}u,
D^{2}u\| _{L_{p}(\bR^{d+1})}
$$
\begin{equation}
                                      \label{12.10.1}
+N R^{ (d+2)(1/p-1/q)-1}\|u\|_{L_{p}(\bR^{d+1})}.
\end{equation}
\begin{theorem}
                                      \label{theorem 12.9.1}
Under the conditions of Lemma \ref{lemma 12.7.1}
 there exists $\hat b>0$, depending only on 
$ d,d_{0}, \uR $, $p,q$, $\bar b_{\uR}$,  $\|b\|_{L_{d+1}(\bR^{d+1})}$,
and
the function  $\bar N(d,d+1,\cdot)$,
such that, if
$\bar b_{R_{0}}\leq\hat b$ is satisfied for an $R_{0}\in (0,\uR]$, 
then for any 
 $u\in C^{\infty}_{0}(\bR^{d+1})$ and $\lambda\geq0$,
$$
\|\partial_{t}u,
D^{2}u\| _{L_{p}(\bR^{d+1})}+\|Du\|_{L_{q}(\bR^{d+1})}
+(\lambda-N) \| u\|_{L_{p}(\bR^{d+1})}+(\lambda-N)
\| u\|_{L_{q}(\bR^{d+1})}
$$
\begin{equation}
                                      \label{12.9.4}
\leq  
N\|L u-\lambda u\|_{L_{q}(\bR^{d+1})}
+N\|L u-\lambda u\|_{L_{p}(\bR^{d+1})},
\end{equation}
where $N$  depend only on 
$ d,d_{0}, \uR $, $p,q$, $R_{0}$, $\bar b_{\uR}$, $\bar b_{\infty}$, 
$\|b\|_{L_{d+1}(\bR^{d+1})}$, and
the function  $\bar N(d,d+1,\cdot)$

\end{theorem}

Proof. By Theorem 5.2 of \cite{Kr_21_1} for $\lambda\geq1$
$$
 \lambda  \| u\|_{L_{p}(\bR^{d+1})}
\leq N\|L u-\lambda u\|_{L_{p}(\bR^{d+1})},\quad
\lambda 
\| u\|_{L_{q}(\bR^{d+1})}
\leq N\|L u-\lambda u\|_{L_{q}(\bR^{d+1})}.
$$
It follows that it suffices to prove \eqref{12.9.4}
for $\lambda=0$.

By classical results
$$
\|\partial_{t}u,
D^{2}u\| _{L_{p}(\bR^{d+1})}\leq N\|L_{0}u
\| _{L_{p}(\bR^{d+1})}\leq N\|L u
\| _{L_{p}(\bR^{d+1})}+N_{2}\|b^{i}D_{i}u\|_{L_{p}(\bR^{d+1})},
$$
where the last term by H\"older's inequality is dominated by
the product
$ 
\|b\|_{L_{d+1}(\bR^{d+1})}
\|Du\|_{L_{q_{p}}(\bR^{d+1})}$ and,
for any $\varepsilon>0$,
$$
\|Du\|_{L_{q_{p}}(\bR^{d+1})}
\leq \varepsilon 
\|Du\|_{L_{p^{*}}(\bR^{d+1})}+
N(\varepsilon)\|Du\|_{L_{q }(\bR^{d+1})} 
$$
\begin{equation}
                                 \label{12.11.5}
\leq N_{3}\varepsilon\|\partial_{t}u,
D^{2}u\| _{L_{p}(\bR^{d+1})}+N(\varepsilon)
\|Du\|_{L_{q }(\bR^{d+1})},
\end{equation}
where $p^{*}=p(d+2)/(d+2-p)$ and the last inequality
is a consequence of   embedding theorems.
Hence,
$$
N_{2}\|b^{i}D_{i}u\|_{L_{p}(\bR^{d+1})}
 \leq N_{3}\varepsilon\|\partial_{t}u,
D^{2}u\| _{L_{p}(\bR^{d+1})}+N(\varepsilon)
\|Du\|_{L_{q }(\bR^{d+1})}.
$$

We also take into account \eqref{12.10.1}
and conclude that
$$
\|\partial_{t}u,
D^{2}u\| _{L_{p}(\bR^{d+1})}\leq  N\|L u
\| _{L_{p}(\bR^{d+1})}+N_{3}\varepsilon\|\partial_{t}u,
D^{2}u\| _{L_{p}(\bR^{d+1})}
$$
$$
+N(\varepsilon)\Big(
R  \| Lu\|_{L_{q}(\bR^{d+1})} 
+  
R^{\alpha}
\|\partial_{t}u,
D^{2}u\| _{L_{p}(\bR^{d+1})}
+  R^{ -\beta}\|u\|_{L_{p}(\bR^{d+1})}\Big),
$$
where $\alpha>1$ and $\beta>0$ are obviously defined 
quantities. We choose and fix $\varepsilon>0$
so that $N_{3}\varepsilon\leq 1/4$ and after that
we make the final choice for  $\hat b$ and
 $R_{0}$ by requiring
not only $N_{1}\bar b_{R_{0}}\leq 1/4$, but also
$N(\varepsilon)R_{0}^{\alpha}\leq 1/4$. Then we get
$$
\|\partial_{t}u,
D^{2}u\| _{L_{p}(\bR^{d+1})}\leq  N\|L u
\| _{L_{p}(\bR^{d+1})} +N\|L u
\| _{L_{q}(\bR^{d+1})} 
+  N|u\|_{L_{p}(\bR^{d+1})}.
$$
After that it only remains to use \eqref{12.10.1}
again. The theorem is proved.

{\bf Proof of Theorem \ref{theorem 2.16.1}}. Uniqueness follows 
from Theorem 5.2 of \cite{Kr_21_1}.
To prove the existence, first assume that $b$ is
bounded and smooth. Then by classical results
(for any $\lambda>0$) we have a unique solution
in $W^{1,2}_{p}(\bR^{d+1})\cap W^{1,2}_{q}(\bR^{d+1})$.
Then take $\zeta\in C^{\infty}_{0}(\bR^{d+1})$
with unit integral and support in the unit ball, 
for $n=1,2,...$ define $\zeta_{n}(z)
=n^{d+1}\zeta(nz)$, and let $b_{n}=b*\zeta_{n}$.
Observe that the quantities $\bar b_{R}$ remain the
same for all $b_{n}$. Therefore, for an appropriate $N_{0}$,
for the solution
$u_{n}$ of $(L_{0}+b^{i}_{n}D_{i})u_{n}-\lambda u_{n}=f$
for $\lambda>N_{0}$ we obtain uniform estimates
of the left-hand sides of \eqref{12.9.4}
with $u_{n}$ in place of $u$. By 
finding $u$ with the properties in \eqref{12.11.4}
and
 a subsequence
$n'$ such that $\partial_{t}u_{n'},D^{2}u_{n'}\to 
\partial_{t}u,D^{2}u $ weakly in $L_{p}(\bR^{d+1})$,
$Du_{n'}\to Du$ weakly in $L_{q}(\bR^{d+1})$
and in $L_{q_{p}}(\bR^{d+1})$ (see \eqref{12.11.5}), and
$u_{n'}\to u$  weakly in $L_{p}(\bR^{d+1})$, and observing
that then $b^{i}_{n'}D_{i}u_{n'}\to b^{i}D_{i}u$
weakly in $L_{p}(\bR^{d+1})$, we easily pass to the limit
in $L_{0}u_{n'}+b^{i}_{n}D_{i}u_{n'}-\lambda u_{n'}=f$.
The theorem is proved.

As a corollary of this theorem we obtain the following
result about solvability of terminal-value
problem with zero data at the finial time.
This result is obtained just by taking $f(t,x)=0$
for $t\geq T$ and multiplying functions by $e^{\lambda t}$.

\begin{theorem}
                               \label{theorem 11.12.3}
Let $T\in(0,\infty)$,
 $p\in [d_{0}+1,d+1)$, 
$q>q_{p}$,
 and let
 $f\in L_{p}((0,T)\times\bR^{d})
\cap L_{q}((0,T)\times\bR^{d})$. 
Assume that condition \eqref{2.18.1} is satisfied. Then
  there exists a unique
solution of the equation $Lu =f$
in $(0,T)\times\bR^{d}$ with terminal condition
$u(T,\cdot)=0$
 in the class of functions
such that
$$
\partial_{t}u,D^{2}u\in L_{p}((0,T)\times\bR^{d}),\quad
Du\in L_{q}((0,T)\times\bR^{d}),
$$
$$
u\in L_{p}((0,T)\times\bR^{d})
\cap L_{q}((0,T)\times\bR^{d}).
$$
\end{theorem}

\begin{remark}
                                   \label{remark 2.18.2}
In Theorem \ref{theorem 11.12.3} the terminal
data is zero. One can easily consider more general
data, say $g(x)$ such that there exists 
$g\in (W^{1,2}_{p}\cap W^{1,2}_{q})((T,T+1)\times
\bR^{d})$ such that $g(T,x)=g(x)$ and $g(T+1,x)=0$.
Indeed, then one would apply Theorem \ref{theorem 11.12.3}
with $T+1$ in place of $T$ to $b(t,x)I_{t<T}$ in place of $b$
and $fI_{t<T}+(\partial_{t}g+\Delta g)I_{t>T}$
in place of $f$.

\end{remark}

\mysection{Application to It\^o's equations}

As we know from \cite{Kr_21_1} there are weak solutions
of the equation
\begin{equation}
                                       \label{12.11.8}
x_{t}=w_{t}+\int_{0}^{t}b(s,x_{s})\,ds,
\end{equation}
where $w_{t}$ is a $d$-dimensional Wiener process.

\begin{theorem}
                                  \label{theorem 12.11.5}
Assume \eqref{2.18.1} and $\|b\|_{L_{d+1}(\bR^{d+1})}<\infty$. Then all solutions
of \eqref{12.11.8} have the same finite-dimensional
distributions.
\end{theorem}

Proof. As is quite common we will rely on solutions
of the corresponding parabolic equations and use It\^o's
formula. The only difficulty is that we have to prove
that the formula is applicable with our $u$ and $b$.

Take $T>0$ smooth bounded $f(t,x)$ with compact support
and let $u$ be the function from Theorem \ref{theorem 11.12.3}.
It is convenient to extend $u(t,x)$ for $t>T$ as zero.
This was actually the way it was meant to be  constructed,
by solving the equation with $f(t,x)=0$ for $t>T$.
Introduce $u_{n}=\zeta_{n}*u$. By It\^o's formula
for any   stopping time $\tau\leq T$ and $t\leq\tau$
$$
u_{n}(\tau,x_{\tau})=
u_{n}(t,x_{t})+\int_{t}^{\tau}Lu_{n}(s,x_{s})\,ds
+\int_{t}^{\tau} D u_{n}(s,x_{s})\,dw_{s}.
$$
As for any stochastic integral one can choose
 $\tau_{n}\uparrow T$ such that
$$
E\Big(\int_{t}^{\tau_{n}} D u_{n}(s,x_{s})\,dw_{s}
\mid \cF^{x}_{t}\Big)=0.
$$
In that case
\begin{equation}
                                         \label{12.12.1}
E\big(u_{n}(\tau_{n},x_{\tau_{n}})\mid \cF^{x}_{t}\big)
=u_{n}(t,x_{t})+E\Big(\int_{t}^{\tau_{n}}Lu_{n}(s,x_{s})\,ds
\mid \cF^{x}_{t}\Big).
\end{equation}

According to Theorem 4.9 of \cite{Kr_21_1}
$$
E \int_{t}^{\tau_{n}}|Lu_{n}-Lu|(s,x_{s})\,ds
\leq E \int_{t}^{T }|Lu_{n}-Lu|(s,x_{s})\,ds
$$
$$
\leq N(1+T)^{\chi}\|\Phi_{1/T}(Lu_{n}
-Lu)\|_{L_{p}((0,T)\times \bR^{d } )},
$$
where   $\Phi _{\lambda}(t,x)=\exp(- \sqrt\lambda 
(|x|+ \sqrt t)\theta )$ and $\chi$ and $\theta$ are independent
of $T$, $u_{n}$, $u$. Observe that
 $\partial_{t}u_{n},D^{2}u_{n}
\to \partial_{t}u ,D^{2}u $ in $L_{p}( (0,T)\times \bR^{d })$
and, owing to the fact that $Du_{n}\to Du$ in any
$L_{r}( (0,T)\times \bR^{d })$ ($f\in L_{r}( \bR^{d+1})$) for any
$r\in [q_{p},\infty)$, we also have (by H\"older's
inequality)
$b^{i}D_{i}u_{n}\to b^{i}D_{i}u$ in
$L_{p}(\Gamma)$, where $\Gamma
\subset  (0,T)\times \bR^{d }$ is any standard
cylinder with base that is unit ball. It follows that
$\Phi_{1/T}(Lu_{n}-Lu)\to 0$ in $L_{p}( (0,T)\times \bR^{d })$.

By similar reasons 
$$
  E \int_{t}^{T }| Lu|(s,x_{s})\,ds
\leq N(1+T)^{\chi}\|\Phi_{1/T} Lu \|_{L_{p}((0,T)\times \bR^{d })}<\infty.
$$
Also observe that by embedding theorems $u$
in $t\leq T$, and, hence,
$u_{n}$ in $t\leq T $ are uniformly continuous
(even H\"older continuous since $p>d/2+1$). This allows us
to pass to the limit in \eqref{12.12.1} and conclude that
\begin{equation}
                                          \label{12.12.2}
-u (t,x_{t})=E\Big(\int_{t}^{T}f(s,x_{s})\,ds
\mid \cF^{x}_{t}\Big).
\end{equation}
Now suppose that for some $n=0,1,2,...$, any
$0\leq t_{0}\leq t_{1}\leq...\leq t_{n}$ (no $t_{1}$
if $n=0$)
and any continuous   $f_{0}(x),...,f_{n}(x)$
with compact support
the quantity
$$
Ef_{0}(x_{t_{0}})\cdot...\cdot f_{n}(x_{t_{n}})
$$
is independent of which solution of \eqref{12.11.8}
we take. Automatically, of course, the same holds if $f_{k}$'s
are just bounded and continuous.

This induction hypothesis holds true for $n=0$,
because \eqref{12.12.2} with $t=0$ implies that
$$
 \int_{0}^{T}E f(s,x_{s})\,ds\quad\text{and, hence,}
\quad E  f(s,x_{s})
$$
is independent of which solution of \eqref{12.11.8}
we take.

To show that the induction works observe that
 by using \eqref{12.12.2}
with $t=t_{n}$, $T=t_{n+1}$ we get that
$$
Ef_{0}(x_{t_{0}})\cdot...\cdot f_{n}(x_{t_{n}})
\int_{t_{n}}^{t_{n+1}}f(s,x_{s})\,ds
$$
is independent of which solution of \eqref{12.11.8}
we take. As above this leads to the conclusion that
$$
Ef_{0}(x_{t_{0}})\cdot...\cdot f_{n}(x_{t_{n}})
 f(t_{n+1},x_{t_{n+1}}) 
$$
is independent of which solution of \eqref{12.11.8}
we take. This, obviously, proves the theorem.

{\bf Acknowledgment}. The author
is sincerely grateful to Hongjie Dong and M. Safonov
for their comments, corrections, and discussions
relevant to the article.

\end{document}